\documentclass[reqno]{amsart} \usepackage{amsfonts}
\usepackage{amssymb} \usepackage[T1]{fontenc} \usepackage{ae}
\usepackage{hyperref}

\numberwithin{equation}{subsection}

\title{The 2-modular permutation modules on fixed point free involutions of symmetric groups}
\author{Peter Collings}
\begin{document}
\maketitle
\begin{center}
  {27 Lillington Road Leamington Spa Warwickshire \newline}
  {pmfcollings@hotmail.co.uk}
\end{center}
\baselineskip=12pt

\newcounter{Subsection}[section] \newcounter{chapterno}
\setcounter{chapterno}{1} \newcounter{lastplace}
\newcommand{\place}{{\bf \arabic{chapterno}.\arabic{lastplace}}}
\newcommand{\Subsection}[2]{\addtocounter{Subsection}{1}
  \setcounter{lastplace}{} \newpage
  #2
\begin{center} \LARGE{\arabic{section}.\arabic{Subsection}}\ \ \ \ \ {\LARGE\texttt {#1}} \end{center}   
\addtocounter{lastplace}{1}}
    \newenvironment{theorm}[2] {
    \noindent \medskip \\ {\large{(\place)}} {\bf { Theorem.}} {\it { #1}}
    \noindent \medskip \\ {\it Proof.} #2$\! \! \! \! \! \qed$}{\noindent 
    \addtocounter{lastplace}{1}}
    \newenvironment{lemma}[2] {
    \noindent \medskip \\ {\large{(\place)}} {\bf { Lemma.}} {\it { #1}}
    \noindent \medskip \\ {\it Proof.} #2$\! \! \! \! \! \qed$}{\noindent 
    \addtocounter{lastplace}{1}}
    \newenvironment{corollary}[2] {
    \noindent \medskip \\ {\large{(\place)}} {\bf { Corollary.}} {\it { #1}}
    \noindent \medskip \\ {\it Proof.} #2$\! \! \! \! \! \qed$}{\noindent 
    \addtocounter{lastplace}{1}}
    \newenvironment{remark}[1] {
    \noindent \medskip \\ {\large{(\place)}} {\bf { Remark.}} {#1}}
  {\noindent \addtocounter{lastplace}{1}}
  \newenvironment{block}[2] {
    \noindent \medskip \\  
    #2{\large{\place}} { #1}} {\noindent \addtocounter{lastplace}{1}}
  
  \newenvironment{fact}[1] {
    \noindent \medskip \\  
    {\large{(\place)}} {\bf { Fact.}}{\it { #1}}} {\noindent
    \addtocounter{lastplace}{1}}
  \newenvironment{definition}[1] {
    \noindent \medskip \\  
    {\large{(\place)}} {\bf { Definition.}} #1} {\noindent
    \addtocounter{lastplace}{1}}
  \newenvironment{convention}[1] {
    \noindent \medskip \\  
    {\large{(\place)}} {\bf { Convention.}} #1} {\noindent
    \addtocounter{lastplace}{1}}
  \newenvironment{example}[1] {
    \noindent \medskip \\ {\large{\place}} {\bf { Example.}} {#1}}
  {\noindent \addtocounter{lastplace}{1}}
  \newenvironment{nblock}[2] {
    \setcounter {lastplace}{1}
     \noindent \medskip \\ \begin {center}
       { \large { \arabic {chapterno}.}} { \large \it{#2} } \newline
     \end {center} { #1} }
   
   \newcommand{\zerostart}{\addtocounter{Subsection}{-1}}
   \newcommand{\zeroblock}{\addtocounter{lastplace}{-1}}
   \setcounter{lastplace}{1}
   \newcommand{\mylabel}[1] { \newcounter{#1}
    \setcounter{#1}{\arabic{lastplace}} \newcounter{#1_cn}
    \setcounter{#1_cn}{\arabic{chapterno} } }
   \newcommand{\partn}{ {\tiny{ | \! - } } }

   \newcommand{\arabicc}[2]{(\arabic{#1_cn}.\arabic{#1})}
   \newcommand{\Fix}{{\rm Fix}} \newcommand{\Stab}{{\rm Stab}}
   \newcommand{\supp}{{\rm supp}} \newcommand{\Sym}{{\rm Sym}}
   \newcommand{\Alt}{{\rm Alt}}
   \newcommand{\alt}{{\rm alt}}
   \newcommand{\ind}{{\rm ind}}
   \newcommand{\Ind}{{\rm Ind}}
   \newcommand{\Syl}{{\rm Syl}} \newcommand{\SM}{{\mathcal S}}
   \newcommand{\SQ}{{\mathcal S ^ q}} \newcommand{\F}{{\mathcal F}}
   \newcommand{\Od}{{\mathcal O}}
   \newcommand{\sgn}{{\rm sgn}}

\markright{$p$-Permutation Modules}
\begin{block}
  { Let $ A = k G $, the group algebra of some finite group where the
    characteristic of the field $ k $ divides $ | G | $. In contrast
    to working over the complex field, the $ k G $-modules are not
    usually semisimple. If a Sylow $ p $-subgroup of $ G $ is not
    cyclic then there are infinitely many indecomposable $ k G
    $-modules, and we usually enjoy little control over the category
    of such modules. It is therefore instructive to find classes of
    modules which may be expressed as a sum of a not very great number
    of indecomposables, and to understand the structure of these
    indecomposables.  Permutation $ k G$-modules, and their
    indecomposable summands (called \linebreak $ p $-permutation
    modules), provide one such class.
        \par
        In \cite{phd04} we studied permutation modules for symmetric
        groups acting on conjugacy classes of fixed point free
        elements which are products of $ q $-cycles, and where $ k $
        has characteristic $ p $. With each component
        (that is, indecomposable summand) we associated a
        fixed point set, and we obtained a general description for the
        fixed point sets. In this paper we shall complete the analysis
        of the case $ q = p = 2 $. Thus we shall determine the fixed
        point sets of the components of the permutation module of the
        symmetric group $ \Sym ( 2 n ) $ on the set of its
        fixed point free involutions.  At the same time we find the
        vertex and Brauer quotient for each component, and 
        the ordinary character associated with each component.
        
        \par 
        \vspace{5mm} Each component
        $ M $ of a permutation
        $ k G $-module determines two
        invariants, namely a vertex $ Q $ (up to conjugation)
        and a projective indecomposable \linebreak $ p $-permutation module for
        for $ N_G ( Q ) / Q $; these invariants characterise $ M $
        up to isomorphism. These invariants were introduced by Green
        \cite{gre58} and adapted to
        permutation modules by Puig; see the paper of
        Brou{\'e} \cite{bro85}.
        \par 
        Given an indecomposable $ k G $-module $ M $, a vertex of $ M
        $ is a subgroup $ H \leq G $ minimal subject to satisfying the
        condition that there be a $ k H $-module $ L $ for which $ M $
        is a component of $ {\rm\Ind} ^ G _ H ( L ) $. This vertex is
        unique up to conjugation in $ G $. Any Sylow $ p $-subgroup of
        $ G $ satisfies this condition and so vertices are always $ p
        $-subgroups.  The $ k H $-module $ L $ is called a source of $
        M $ and is unique up to conjugation in $ N_G ( H ) $. However,
        $ M $ is a $ p $-permutation module if and only if $ M $ is a
        component of a permutation module and this happens exactly
        when $ L $ is the trivial $ k H $-module. Hence $ p
        $-permutation modules are also known as trivial source modules.

        \par 
        The Brauer quotient $ M ( H ) $ of $ M $ is a
        module for the group $ N_G ( H ) / H $. When $ M $ is a
        component of a permutation module $ k \Omega $ then $ M ( H ) $
        is a component of $ k {\rm Fix}_{\Omega} H $, hence is a $ p
        $-permutation module. Brou\'e's correspondence (1.1) asserts a
        $ 1 $-$ 1 $ correspondence between the projective components of $
        k {\rm Fix}_{\Omega} H $ and the components of $ M $ with
        vertex $ H $.

        \par
        In \cite{phd04} we classified the possible fixed
        point sets when the group is a symmetric group acting by
        conjugation. We will summarize the relevant results
        from \cite{phd04} in
        \S 1 below.  We are then in a position to determine the
        fixed point sets, and the vertices of the components of, the
        permutation module afforded by the action of $ \Sym ( 2 n ) $
        on its conjugacy class $ \Xi_{ 2 n } $ of fixed point free
        involutions. This permutation action has been studied before;
        as an example the ordinary character it affords is known, and
        we refer to \cite{ing90} for a demonstration of this and other
        results. With any $p$-permutation module one can
        associate an ordinary character, and we will determine these characters
        for the components of the permutation module on  fixed point free 
        involutions.
        For general results we refer to \S 26 and \S 27 in \cite{the95}.
        }  { }
\end {block}
\setcounter{chapterno}{1} \setcounter{lastplace}{1}
\vspace*{5mm}
\begin{center} {\large {\bf \arabic{chapterno}.} {\bf Summary of results from \cite{tp05}} } \end{center}
\vspace*{-7mm}
\begin{block}
{
        {\bf Fixed point sets} Let $ G $ be an arbitrary
        finite group, $ \Omega $ some finite $ G $-set and $ M = k \Omega
        $ the resulting permutation module. In \cite{bro85} the
        components of $ M
        $ which have a vertex $ Q $ are parameterized as follows.}
{}
\end{block}        
\mylabel{BCorr}
\begin{block}
  {  {\bf Theorem } {\it Let $ \Omega $ be a permutation $ G $-space. There is a
      multiplicity-preserving 1-1 correspondence between \newline (i)
      components of $ M $ with vertex $ Q $, and \newline (ii)
      components of $ M ( Q ) $ which are projective as modules for $
      N_G ( Q ) / Q $. }  } {}
\end{block}

Here $ M ( Q ) = k [ \Fix ( Q ) ] $, which is a permutation module for
the group $ N_G ( Q ) / Q $ in the natural way, and $ \Fix _ \Omega (
Q ) $ is the set of fixed points $ Q $ in $ \Omega $. We will refer to
this as `Brou\' e Correspondence'.

\par
The main object of study in \cite{tp05} are the subsets of $ \Omega $
which can occur as fixed point sets. The basic definition is as
follows.
\mylabel{definition:fps}
\begin {definition}
  { {\it Suppose $ X $ is a subset of $ \Omega $. Then $ X $ is a
      fixed point set if $ X = \Fix_\Omega ( Q ) $ where $ Q $ is a
      vertex for some component $ W $ of $ k \Omega $.}  }
\end {definition}

For any subset $ X $ of $ \Omega $, we write $ S_X := \Stab_G ( X ) $
for the pointwise stabilizer of $ X $, and we write $ Q_X $ for some
Sylow $ p $-subgroup of $ S_X $.  \par We say that the subset $ X
$ of $ \Omega $ is {\bf closed} if $ \Fix_\Omega ( Q_X ) = X $.
\par It is easy to see that any fixed point set is closed,
furthermore the group $ Q $ in the above definition of a fixed point
set can be taken as $ Q = Q_X $.  [See \cite{tp05} (2.3)] Now let $
N_X $ be the set stabilizer of $ X $, and let $ {\overline N_X} := N_X / S_X $;
then $ X $ is a permutation $ {\overline N_X} $-set. Then we have the following
modification of the Brou\' e Correspondence (but not yet taking
multiplicities into account).  [See \cite{tp05} (2.12)]

\begin{block}
  { {\bf
      Proposition.} Assume $ \Omega $ is a $ G $-set, and $ X \subseteq \Omega $ is a
    closed subset. Then $ X $ is a fixed point set if and only if $ k
    X $ as a module for $ {\overline N_X} $ has a projective component.  }  {   }
\end{block}

\begin{block}
  {  {\bf Fixed point sets by conjugation} We consider now a symmetric group acting by conjugation
    on itself. In this case, $ G $-sets derive additional structure
    from multiplying permutations. The general setup studied in
    \cite{tp05} is as follows.  \par We consider permutations with
    finite supports on arbitrary subsets of the natural numbers, that
    is, let $ \Gamma $ be the group of finitary permutations on $
    {\mathbb N} $. Let
    $$
    {\mathcal F} := \{ X \subset \Gamma , | X | < \infty \}
    \setminus \{ id \}.
    $$
    Most important is the subset $ \SQ \subset \F $, where $ X $ is
    an element of $ \SQ $ if and only if it consists of permutations
    which are products of $ q $-cycles.  In our application later we
    will take $ q = 2 $.  \vspace*{1mm} \par For $ X \in \F $, let
    $ G_X $ be the group of all permutations on the support $ \supp (
    X ) $. Then we have a natural $ G_X $-set which has $ X $ as a
    subset, namely
    $$
    \Xi_X := \cup_{ g \in G_X } { X ^ g } .$$
    We want to understand
    when $ X $ is a fixed point set, as a subset of $ G_X
    $-permutation set $ \Xi_X $.    
  } { }
\end{block}
\begin{block}
  { {\bf Multiplicative structure of $ \F $} We define an equivalence relation on $ \F $.  Take $ X_1
    $ and $ X_2 $ in $ \F $. We say
    $$
    X_1 \sim X_ 2 \Leftrightarrow {\rm \ there \ is} g \in \Gamma
    {\rm \ such \ that \ } X_1 ^ g = X_2.
    $$
    \par We work with equivalence classes but keep the same notation.
    (See \cite{phd04} for details.)  \vspace*{1mm} \newline {\bf
      Constructions} (1) [Products] Let $ X , Y \in \F $. Take $ g \in
    \Gamma $ such that $ X $ and $ Y ^ g $ have disjoint supports.
    Define
    $$
    X * Y := {\rm \ the \ class \ containing \ } X * Y ^ g .$$
    \newline (2)
    [Powers] Let $ X \in \F $.  Define $ * ^ 2 X = X * X $ and
    inductively
    $$
    * ^ s X = ( * ^ { s- 1 } X ) * X $$
    which is the class of $ X
    \times X ^ { g_2 } \times \dots \times X ^ { g_s } $ for $ g_i \in
    \Gamma $ where the fators have disjoint supports. We write $ X ^ s
    $ for simplicity.  \newline (3) [Diagonals] Let $ X \in \F $.
    Then define $ \Delta ^ s ( X ) $ to be the class of a subset of $
    * ^ s X $. If $ * ^ s X $ is represented by $ X \times X ^ { g_2 }
    \times \dots \times X ^ { g_s } $ as in (2), then $ \Delta ^ s ( X
    ) $ is the class of the subset
    $$
    \{ x * x ^ { g_2} * \dots * x ^ { g_s } | x \in X \} . $$
    \par We say that an element $ X \in \F $ is {\bf irreducible}
    if there are no $ Y , Z \in \F $ such that $ X = Y * Z $. Then
    every $ X \in \F $ has a factorisation into a product of
    irreducible sets, and this is essentially unique. Given two
    elements $ X $ and $ Y $ in $ \F $, we may say that $ X $ and $ Y
    $ are {\bf coprime} if they do not share an irreducible factor
    (upto equivalence).  \par A set $ X \in \F $ is called {\bf
      exact} if it is fixed point free on the support $ \supp ( X ) $.
    On the other hand we call an irreducible element $ X $ 
    {\bf projective} if $ Q_X = 1 $, and an element $ Y $
    projective-free if none of its irreducible factors is projective.
    It is proved in \cite{phd04} that if $ X $ is an irreducible fixed point set
    then it is either exact or projective. The second case
    occurs precisely when $ X $ corresponds to a projective component of
    the permutation module.  \par Furthermore, a fixed point set
    is said to be {\bf transitive} if the group $ \langle X \rangle $
    is transitive on the support of $ X $. \par Finally,
    for $ X \in \F $ let the number $ \kappa ( X ) $
    be the lowest positive
    integer $ u $ such that the permutation module $ k X ^ {\wr u} $
    for $ H \wr \Sym ( u ) $ does not have a projective component, allowing
    the possibility $ \kappa ( V ) = \infty $ if no such $ u $ exists.
    \par With these, the
    general description of a fixed point set in $ \SQ $ given in
    \cite{phd04} [Theorem 7.18] is as follows.  Let $ Y \in \SQ $,
    write $ Y = Y_1 ^ { a_1 } * Y_2 ^ { a_2 } * \dots * Y_t ^ { a_t }
    $ where the $ Y_i $ are pairwise coprime and irreducible.  }  {
     }
\end{block}
\mylabel{theorem:mainphd}
\begin{block}
  { {\bf Theorem} [[3], 7.18] {\it Assume the setup of before, with elements of $
      \mathcal S ^ q $.  Assume $ X \in \mathcal S ^ q $.  \newline
      (1) $ X $ is an irreducible exact fixed point set $
      \Leftrightarrow $ $ X = \Delta ^ { p^i } Y $ where $ Y $ is a
      transitive irreducible exact fixed point set and $ i \geq 0 $.
      \newline (2) $ X $ is a projective-free fixed point set $
      \Leftrightarrow $ $ X = X_1 ^ { a_1 } * \dots * X_t ^ { a_t } $
      where the $ X_i $ are pairwise coprime irreducible exact fixed
      point set, and $ 1 \leq a_i < \kappa ( X_i ) $.  \newline (3) $
      X $ is a fixed point set $ \Leftrightarrow $ $ X = W * V $ where
      $ W $ is a projective-free fixed point set and $ V $ is an
      irreducible projective fixed point set.  }  }  {  }
\end{block}
The motivation for considering transitivity is the following (see
\cite{phd04}, 7.14). For $ Y \in \SQ $ we define the {\bf degree} of $
Y $ to be $ {\rm d} ( Y ) = | \supp ( Y ) | $.
\mylabel{theorem:transitive}
\begin{block}
  { {\bf Theorem} [[3], 7.14] {\it Assume $ Y \in \SQ $ is an exact fixed point set
      which is transitive. Then $ {\rm d} ( Y ) = p $ or $ pq $.  }  }
  {  }
\end{block}
\vspace*{0.1cm}
\par This will be sufficient to classify all fixed point sets which occur
for components of our fixed point free permutation module in the case $
p = q = 2 $.

\setcounter{chapterno}{2}
\setcounter{lastplace}{1}
\vspace*{5mm}
\begin{center} {\large {\bf \arabic{chapterno}.} {\bf The components of $ k \Xi_{ 2 n } $} } \end{center}
\vspace*{-7mm}
\begin{block}
  { Let $ k $ be a field of characteristic $ 2 $ and $ \Xi_{ 2 n }
    $ be the conjugacy class of $ \Sym ( 2 n ) $ containing the fixed
    point free involutions. We intend to enumerate the components of
    the permutation $ k ( \Sym ( 2 n ) )$-module $ k \Xi_{ 2 n }
    $, along with the associated vertices and Brauer characters, by
    enumerating its fixed point sets.  As we have mentioned this
    permutation action has been studied before. The
    ordinary character afforded by this action has been calculated in
    general and the paper \cite{ing90} contains a proof of this
    result. We quote the following theorem from the same paper which
    provides an even more general result, which will be useful later
    when we determine the characters of the components of $ k \Xi_{ 2
      n } $.
        \par 
        To describe it, we let $ n $ and $ m $ be non-negative
        integers.  Let $ \sgn_{2 m} $ be the module over $ {\mathbb C} $ afforded by
        the sign representation of $ \Sym ( 2 m ) $.  Then the outer
        tensor product $ {\mathbb C} \Xi_{2 n} \# \sgn_{2m} $ is a $ \Sym (
        2 n ) \times \Sym ( 2 m ) $-module, and the group $ \Sym ( 2 n
        ) \times \Sym ( 2 m ) $ is a Young subgroup of $ \Sym ( 2
        \cdot ( m + n ) ) $.  Let $ \Xi_{ 2 n , 2 m } := \Ind ^ {
          \Sym ( { 2 ( n + m )} ) }  _ { \Sym ( {2 n } ) \times \Sym ( { 2 m} ) } 
        ( {\mathbb C} \Xi_{2 n} \# \sgn_{ 2 m }) $ denote the induced $ {\mathbb
          C} \Sym ( { 2 ( n + m ) } ) $-module. If $ n = 0 $ then put $
        \Xi_{ 2n , 2m } := \sgn_{ 2 m } $.
        \par
        Note that we have $ \Xi_{ 2n , 0} = \Xi_{ 2n } $, which is
        the permutation $ \Sym ( 2 n ) $-space under consideration.
        
      } {  }
\end{block}
\mylabel{theorem:irs}
\begin {theorm}
  { {\rm Inglis-Richardson-Saxl \cite{ing90}} Let $ \Lambda_{ 2 n } ^
    m $ be the set of partitions of $ 2 n $ using precisely $ m $ odd
    parts.  The character $ \chi_{ 2n , 2m } $ of the $ {\mathbb C}
    \Sym ( { 2 ( n + m ) } ) $-module $ \Xi_{ ( 2n , 2m ) } $ satisfies $
    \chi_{ ( 2n ,2m ) } = \sum_{ \lambda \in \Lambda_{ 2 ( n + m ) } ^ { 2
        m} } \chi ^ { \lambda } $.  }  { A short proof is given in
    \cite{ing90}.  }
\end {theorm}
\mylabel{lemma:irs}
\begin {lemma}
  { {\it For every positive integer $ n $ the complex character $
      \chi_n $ afforded by the permutation $ \Sym ( 2 n ) $-space $
      \Xi_{2 n} $ is given by
      $$
      \chi_n := \chi_{ \Xi_{2 n} } = \sum_{ \lambda \in
        \Lambda_{ 2 n } ^ 0 } \chi ^ { \lambda } .$$} } { This is
    the previous theorem applied to the case $ m = 0 $.  }
\end{lemma}
\mylabel{theorem:murray}
\begin {theorm}
  { The permutation $ k \Sym ( 2 m ) $-module $ k \Xi_{2m} $ admits
    no projective components.  }  { This follows from the main
    theorem of \cite{mur06} (J. Murray). }
\end {theorm}
\mylabel{lemma:natural}
\par            The following lemma presents an elementary result which
will help later to describe the fixed point sets of  $ k ( \Xi_{2 n} ) $.
\begin {lemma}
  { {\bf The natural} $ k \Sym ( 3 ) ${\bf -module}.  Let the field $
    k $ be of characteristic $ p = 2 $ with $ G = \Sym ( 3 ) $. Let $
    P = < ( 1 , 2 ) > $ be the Sylow $ 2 $-subgroup of $ G $ generated
    by the transposition $ ( 1, 2 ) $. There is a natural action of $
    G $ on the set $ X := \{ 1 , 2 , 3 \} $ which affords a
    permutation $ k G $-module $ k X $, called the natural $ k \Sym (
    3 ) $-module.  This module admits a projective component.  }  {
    Since $ | X | = 3 $ is coprime to $ 2 $ we have $ k X = 1_{ k G }
    \oplus E $ is the direct sum of the trivial $ k G $-module and the
    submodule $ E $ spanned over $ k $ by the set $ Y := \{ 1 + 3 , 2
    + 3 \} $. The submodule $ E $ is termed the standard $ k \Sym ( 3
    ) $-module. The dimension of $ E $ is $ 2 $. Note that
    $ P $ acts regularly on $ Y
    $ and this ensures that $ E $ is a projective $ k P $-module,
    from which it follows that
    $ E $ is a projective $ k G $-module.}
\end {lemma}
\newcommand{\SetS}{ {\mathcal SetS } }
\mylabel{block:building}
\begin{block}
{   In light of the fact that \arabicc{theorem:murray}{}
    shows that there are no projective components of $ \Xi_{ 2 n } $,
    the problem of finding the fixed point sets reduces by
    \arabicc{theorem:mainphd}{} to finding the transitive fixed point
    sets.
    \par
    By \arabicc{theorem:transitive}{}, any transitive
    fixed point set has support of
    size $ 2 $ or $ 4 $. We can now find these sets.
    \par  
    {\it The singleton set} \ \ Let $ U := \{ (1 2) \} $.
    This is a fixed point set: the group $ G_U $ is
    the symmetric group of degree $ 2 $ and
    $ U = \Xi_2 $. We have $ S_U = G_U = Q_U $ and
    $ k U $ is the trivial module which has vertex $ Q_U $.
    This is the only way to have a fixed point set with
    support of size $ 2 $.
    \par
    {\it Transitive irreducible fixed point sets supported on $ 4 $
      elements} \ \ We must analyze the permutation $ \Sym ( 4 )
    $ module afforded by the permutation space $ \Xi_4 = \{ ( 1 2 ) (
    3 4 ), ( 1 3 ) ( 2 4 ), ( 1 4 ) ( 2 3 ) \} $.  The permutation module
    is isomorphic to $ \Ind ^ G _ D ( k ) $ where $ D $ is the
      stabilizer of $ ( 1 2 ) ( 3 4 ) $ (say) in $ G = \Sym ( 4 ) $.
      So $ D $ is dihedral of order $ 8 $.  This group $ D $ contains
      the Klein $ 4 $-group $ H = {\rm id} \cup \Xi_4 $ which is
        normal in $ G $, and $ H $ acts trivially on $ \Xi_4 $. Viewed
        as a module for the factor group $ G / H \cong \Sym ( 3 ) $,
        the module $ k \Xi_4 $ is the natural $ 3 $-dimensional
        permutation module. As explained in \arabicc{lemma:natural}{}
        we have $ k \Xi_4 = k \oplus E $ where $ E $ is projective as
        a module for $ \Sym ( 3 ) $. The trivial summnad has vertex $
        D $, and $ D $ has fixed points just the singleton set $ \{ (1
        2 ) ( 3 4 ) \} $. This is equivalent to $ U ^ 2 $ which is not
        transitive or irreducible.  \par The module $ E $ has
        vertex $ H $ as a module for $ G $, and $ \Fix ( H ) = \Xi _ 4
        $. So $ \Xi_4 $ is a fixed point set, and it is irreducible
        and transitive.  
}{}
\end{block}
\begin {theorm}
  { {\bf Fixed Point Sets of the Family $ k ( \Xi_{ 2 n } ) $.}  Write $
    U := \Xi_2 $, $ V:= \Xi_4 $ and $ V_i := \Delta ^ { p ^ i }
    v $. Then a complete list without repetitions of irreducible fixed
    point sets is $ U , V_0 , V_1 , V_2 , \dots $. So if $ X \subset
    \Xi_{ 2 n } $ for some value of $ n $ then $ X $ is an irreducible
    fixed point set if and only if $ X = U $ or $ X = V_i $ for some
    positive integer $ i $.  }
{       By \arabicc{theorem:transitive}{} any transitive fixed
        point set has degree $ 2 $ or $ 4 $.
        By \arabicc{block:building}{} then the only
        transitive fixed point sets are $ \Xi_2 $
        or $ \Xi_4 $.
        Let us put $ U := \Xi_2 $ and
        $ V:= \Xi_4 $.  Now let $ X $ be any irreducible fixed
        point set. By \arabicc{theorem:mainphd}{}
        the set $ X $ satisfies $ X = \Delta ^ {
          p ^ i } U $ or $ X = \Delta ^ { p ^ i } V $, for some
        positive integer $ i $. But $ | U | = 1 $ so $ \Delta ^ { p ^
          i } U = * ^ { p ^ i } U $ is not irreducible. Hence either $
        X = U $ or $ X = \Delta ^ { p ^ i } V $ for some $ i $. Let $
        V_i := \Delta ^ { p ^ i } V $. Then $ V_i $ {\it is} an
        irreducible fixed point set by \arabicc{theorem:mainphd}{}.
        Let us write $
        V_0 := V $. Then a complete list of the irreducible fixed
        point sets is $ U , V_0 , V_1 , V_2 , \dots $.  }  { }
\end {theorm}
\mylabel{theorm:MainTheorem}
\begin {theorm}
  { Let $ U := \Xi_2 $ and let $ V_i := \Delta ^ { 2 ^ i } \Xi_4
    $. Let $ I $ be a finite set of non-negative integers and $ s $ a
    non-negative integer. Then $ U ^ s * ( *_{ i \in I } V_i ) $ is a
    fixed point set, assuming that either $ s \neq 0 $ or $ I \neq
    \emptyset $. Furthermore any fixed point set has this form.  }  {
    First we observe that the singleton set $ U ^ s := ( \Xi_2) ^ s =
    \{ ( 1 \ 2 ) ( 3 \ 4 ) \cdots ( 2 s - 1 \ 2 s ) \} $ is always a
    fixed point set and corresponds to the trivial component of $ k
    \Xi_{ 2 s } $, which can be seen directly.
            \par The sets $ V_i $ are pairwise distinct irreducible
            fixed points, hence pairwise coprime,
            and each is coprime to $ U ^ s $.
        \par It follows from \arabicc{theorem:mainphd}{}
        that $ U ^ s * ( *_{ i \in I } V_i ) $ {\it is} a fixed point set.
        Furthermore every fixed point set is
        of the form $ U ^ s * ( *_{ i \in I }
        V_i ^ { a_i }) $ for non-negative integers $ s $ and $ a_i \in
        I$ for some set $ I $ satisfying $ a_i < \kappa ( V_i ) $.  It is not hard to show
        that $ \kappa ( V_i ) = 1 $ for each $ i $. }
\end {theorm}
\mylabel{lemma:permiso}
\begin {lemma}
{
    Let $ X, Y, Z \in \F $ and assume $ X $ irreducible. Then
    the pair $ ( {\overline N}_{ \Delta ^ t X }, \Delta ^ t X ) $
    is permutation isomorphic
    to the pair $ ( {\overline N}_X , X ) $. Furthermore if $ Y $ and $ Z $
    are coprime then $ N_ { Y * Z } = N_{ Y } \times N_{ Z } $.
}
{
    See \cite{tp05} (4.9) and (5.4).
}
\end {lemma}
\mylabel{block:Brauer}
\begin {block}
  { {\bf Brauer Quotients and the Brou\'e Correspondence.} We have determined the irreducible fixed
    point sets of $ k (\Xi_{2 m} ) $. Every such set $ W $
    determines a vertex $ Q_W $ and gives rise to a Brauer quotient.
    We recall that
    this Brauer quotient is the pair consisting of the factor group $ N_{\Sym(2m)}
    (Q_W) / Q_W $ and the permutation module $ k W $ on which the
    factor group acts. The Brou\'e correspondence \arabicc{BCorr}{} tells us that the
    set of projective components of $ k W $ is in $1$-$1$ correspondence
    (preserving multiplicities) with the components of $ k (\Xi_{2 m} ) $ with vertex $ Q_W $. We will
    later show that these sets are singletons.
    The group $ {\overline N_W} := N_W / S_W $ also acts on the set $ W $, and in fact the
    projective components of the module afforded are in $1$-$1$
    correspondence with the projective components of the Brauer quotient.
    This follows because 
    the group $ {\overline N_W} $ is isomorphic
    to the quotient $ (N_{\Sym(2m)} (Q_W) / Q_W) / T $ where
    $ T = N_{ S_W } ( Q_W ) / Q_W $ is a normal
    $ 2 ^ \prime $ subgroup which acts trivially on $ W $. The components of
    $ k W $ as a $ k {\overline N_W} $ module therefore inflate to give the components
    of $ k W $ as a $ N_{\Sym(2m)} (Q_W) / Q_W $ module, and since
    $ T $ is $ 2 ^ \prime $ this preserves projectivity. For this reason we
    will refer to the pair $ ( {\overline N_W} , k W ) $ as the equivalent of the
    Brauer quotient of $ k (\Xi_{2 m} ) $ with respect to $ Q_W $.  
    See \cite{tp05} (2.11) for more details.
    \par
    Consequently it will be sufficient
    for us to describe the permutation module given by the pair $ (
    {\overline N_W}, W ) $, in the sense that the Brou\'e correspondence
    extends to a $1$-$1$ multiplicity preserving correspondence
    between the set of components of $ k (\Xi_{2 m} ) $ with vertex $ Q $
    and the set of projective components of the $ k {\overline N_W} $ module $ k W $.    
    We do this now for the irreducible fixed points sets,
    giving the general result in the theorem that follows.
        \par 
        Firstly consider $ U = \{ ( 1 \ 2 ) \} $ and $ {\overline N}_U = \Sym ( 1 \ 
        2 ) / \Sym ( 1 \ 2 ) \cong {\rm id} $. Here $ {\overline N}_U $ is the trivial
        group and $ k U $ is the trivial $ k {\overline N}_U $-module. In other
        words the pair $ ( {\overline N}_U , U ) $ is permutation isomorphic to
        the pair $ ( 1 , 1 ) $.
        \par 
        Now let $ i $ be a non-negative integer and consider the
        permutation $ {\overline N}_{ V_i } $-space $ V_i $. By
        \arabicc{lemma:permiso}{}         
        pair $ ( {\overline N}_{ V_i } , V_i ) $ is permutation isomorphic to the
        pair $ ( {\overline N}_{ V } , V ) $. Now $ V $ is the set of double
        transpositions in $ \Sym ( 4 ) $ so $ N_V := N_{ \Sym ( 4 ) }
        ( V ) = \Sym ( 4 ) $. Furthermore $ S_V := \Stab_{ \Sym ( 4 ) }
        ( V ) = H $, where $ H $ denotes the Klein $ 4 $-subgroup.
        Thus $ {\overline N}_V \cong \Sym ( 3 ) $ and we are in the situation of
        \arabicc{lemma:natural}{}. In particular the pair $ ( {\overline N}_V , V ) $ is
        permutation isomorphic to the pair $ ( \Sym ( 3 ) , X ) $,
        where $ X $ denotes the natural permutation $ \Sym ( 3 )
        $-space. We saw in \arabicc{lemma:natural}{} that $ k N = k \oplus k E $, the
        summand $ k U $ being indecomposable projective.  }  { }
\end {block}
\markright{Main Results}
\mylabel{theorem:Brauer}
\begin {theorm}
  { Let $ I $ be a finite set of non-negative integers with $ t := | I
    | $.  Let $ s $ be a non-negative integer and put $ W := U ^ s * (
    *_{ i \in I } V_i ) $. Then $ {\overline N_W} $ is
    isomorphic to the direct product of $ t $ copies of $ \Sym ( 3 )
    $. After this identification the permutation $ k {\overline N_W} $-module $ k
    W $ satisfies $ k W \cong \bigotimes ^ t k X $, where $ X $ is the
    natural permutation $ \Sym ( 3 ) $-space.  }  { The set $ U ^ s $
    is a singleton set and so we have $ N_{ U ^ s } = S_{ U ^ s } $,
    and $ {\overline N}_{ U ^ s } = 1 $. Thus the pair $ ( {\overline N}_{ U ^ s } , U ^ s ) $
    is permutation isomorphic to the pair $ ( 1 , 1 ) $.
    By \arabicc{lemma:permiso}{} 
    for any $ i \in I $ the pair $ ( {\overline N}_{ V_i } , V_i ) $ is
    permutation isomorphic to the pair $ ( \Sym ( 3 ) , X ) $.
    The pairwise coprimeness of the irreducible
    factors of $ W $ shows that
    $ N_W = N_{ U ^ s } \times \prod_{ i \in I } N_{ V_i } $, from which
    the result follows.
}        
\end {theorm}
\vspace*{1mm}
\par The aim now is to find the corresponding description of the
vertices.
\begin {convention}
  { Whenever $ W $ is a fixed point set, assume that $ W $ has been
    chosen so that $ W \subset \Xi_{2n} $, where $ n := d ( W ) / 2
    $, so that the support of $ W $ is the set $ \{ 1, 2, \dots , 2 n
    \} $.  This means the subgroup $ Q_W $ is a vertex of some
    component of $ k ( \Xi_{2n} ) $, by \arabicc{definition:fps}{}.  }
  { }
\end {convention}
\mylabel{notation:vertex}
\begin {block}
  { {\bf Notation.}  Fix a positive integer $ n $ and consider the
    permutation $ k \Sym ( 2 n ) $-module $ k ( \Xi_{2 n} ) $. Let $
    \mu := ( s, t ) $ be a pair of non-negative integers satisfying $
    s +2 t = n $. If $ t= 0 $ put $ I := \emptyset $.  Otherwise let $
    t = \sum_{ l = 1} ^ r 2 ^ { i_l } $ be the $ 2 $-adic expansion of
    $ t $ and let $ I := \{ i_l \}_{ l =1 } ^ r $ be the set of
    exponents used.  Define $ W_\mu := U ^ s * ( *_{i \in I} V_i ) $
    and $ Q_\mu := Q_{ W_\mu } $,  $ N_\mu := N_{ W_\mu } $
    and $ { \overline N_\mu } := {\overline N_{ W_\mu } }  $. Then with the above convention in
    force the $ 2 $-subgroup $ Q_\mu $ is a vertex of some
    component of $ k ( \Xi_{2 n} ) $.  }  { }
\end {block}
\vspace*{0.2cm}
\par The size of the support of $ W_\mu $ is given by
$ 2 n = | \supp ( W_\mu ) | = 4 t + 2 s $ where $ 4 t = \sum | \supp ( V_i ) | $.
Since each $ V_i $ has support of size a power of $ 2 $, the possible
$ W_\mu $ are (by uniqueness of $ 2 $-adic expansion)
in $ 1 $-$ 1 $ correspondence with solutions in non-negative integers of the equation
$ 4 t + 2 s = 2 n $.
\mylabel{theorem:summands}
\begin {theorm}
  { Fix a positive integer $ n $ and consider the permutation $ k \Sym
    ( 2 n ) $-module $ k ( \Xi_{ 2 n } ) $. Define the set $ {\mathcal
      Q} := \{ Q_\mu \} $ over all the $ \lfloor n \ 2 \rfloor + 1 $ pairs $ \mu =
    ( s , t ) $ that are non-negative integer solutions to $ 4 t + 2 s =
    2 n $.  Then $ {\mathcal Q} $ is \linebreak \vspace*{0.1cm} a set of
    $ \lfloor n / 2 \rfloor + 1 $ pairwise non-isomorphic $ 2 $-subgroups of $
    \Sym ( 2 n ) $. Furthermore every element $ Q \in {\mathcal Q} $
    is a vertex of exactly one component of $ k ( \Xi_{2 n} ^
    2 ) $.  \linebreak In particular the number of
    such components is $ \lfloor { n } / { 2 }  \rfloor  + 1 $.  }  { By \arabicc{theorem:mainphd}{} the
    general form of a fixed point set $ W = U ^ e * ( * _{j
      \in J} V_j ) $. Here $ e $ is a non-negative integer and $ J $
    is a set of non-negative integers, \linebreak \vspace*{0.1cm} and
    either $ e > 0 $ or $ J \neq \emptyset $. The degree of $ W $ is
    given by
    $$ d ( W ) = | supp \ ( U  )| \cdot e + \sum_{ j \in J} | \supp
    \ ( V_i ) | = 2 e + \sum_{ j \in J} 2 ^ j \cdot 4 .$$ With our
    convention in force, the fixed point set $ W $ satisfies $ W
    \subset \Xi_{ 2 n } $ if and only if $ d ( W ) = 2 n $, and this is
    the case exactly when $ n = e + 2 \cdot ( \sum_{ j \in J} 2 ^ j )
    $. It follows that the set $ {\mathcal Q} $ of vertices of
    components of $ k ( \Xi_{2 n} ) $ is complete.
        \par \vspace*{0.1cm} Suppose that the pairs $ \mu = ( 4 t , 2 s ) $ and $ \lambda = ( 4 t_1 , 2 s_1 ) $ of non-negative integers
        satisfy $ 2 s + 4 t = 2 s_1 + 4 t_1 = 2 n $. Suppose also that $ \mu
        $ and $ \lambda $ are unequal.
        Then $ s $ is the number of orbits on $ \supp ( W_\mu ) $ which
        are of size $ 2 $ under the action of $ \langle W_\mu \rangle $.
        The corresponding remark holds for $ s_1 $ and so
        $ W_\mu $ and $ W_\lambda $ cannot be conjugate in $ \Sym ( 2 n )  $.
        But the sets $ W_\mu $ and $ W_\lambda $ are closed so they
        satisfy $ W_\mu = \Fix \ Q_\mu $ and $ W_\lambda = \Fix \ 
        Q_\lambda $. Thus $ Q_\mu $ and $
        Q_\lambda $ are non-conjugate subgroups of $ \Sym ( 2 n ) $.
        \par \vspace*{0.1cm}
        To see that the subgroup $ Q_\mu $ is a vertex of no more
        than one component of $ ( k \Xi_{ 2 n } ) $ we
        consider the permutation $ k {\overline N}_{ W_\mu } $-module $ k W_\mu $.
        By \arabicc{theorem:Brauer}{} the pair $ ( {\overline N}_{ W_\mu } , W_\mu ) $ is
        permutation isomorphic to the pair $ ( \Sym ( 3 ) ^ t ,
        \bigotimes ^ t k X ) $. By \arabicc{lemma:natural}{} the permutation $ k
        \Sym ( 3 ) $-module $ k X $ admits the unique component
        projective $ k E $, where $ U $ is the standard
        permutation $ \Sym ( 3 ) $-space. A standard result now
        implies that $ \bigotimes ^ t k X $ admits the unique
        projective component $ \bigotimes ^ t k E $. So
        by the Brou\'e correspondence as given in \arabicc{block:Brauer}{}
        $ Q_\mu $ is a vertex of exactly one component.  }
\end {theorm}
\mylabel{block:vertices}
\begin{block}
{
        {\bf Vertices} For each composition $ \mu $ of $ 2 n $
        of the form $ \mu = (4 t , 2 s ) $ we have the fixed point
        set $ W_\mu = ( * _ { i \in I } V_i ) * U ^ s $ where
        $ t = \sum _ { i \in I } 2 ^ i $ is the $ 2 $-adic
        expansion.
        Then by the results of \cite{phd04} a component which
        corresponds to such a fixed point set has vertex
        a Sylow $ 2 $-subgroup of the pointwise stabilizer
        $ S_{ W_\mu } $ of $ W_\mu $. Since the factors
        of $ W_\mu $ are coprime, this stabilizer is just
        the product of the stabilizers of the factors (by \arabicc{lemma:permiso}{}) so that
        $$ S_{ W_\mu } = S _ { U ^ s } \times ( \prod _ { i \in I } S _ { V_i } ) $$
        In the following denote by $ P ( b ) $ a Sylow $ 2 $-subgroup
        of $ \Sym ( b ) $.
        It follows that
        \newline (a) $ Q_ { U ^ s } $ is a Sylow $ 2 $-subgroup $ P ( 2 s ) $
                     of $ \Sym ( 2 s ) $
                     (since $ S _ { U ^ s } $ is the centralizer of a fixed point free involution
                     in $ \Sym ( 2 s ) $.)
        \vspace*{0.1cm}
        \par
        We have (see \cite{phd04}, Lemma (4.5)) that
        $ S _ { V_i } \cong S_V \wr \Sym ( 2 s ) $. Now $ S_V = H $,
        the Klein $ 4 $-group on four points, $ \Xi_4 \subset H $, and
        hence
        \newline (b) $ Q _ { V_i } = H \wr P ( 2 ^ i ) $, on
        $ 4 \cdot 2 ^ i $ points.
        Combining these, a vertex of a component with
        fixed point set $ W_\mu $ where $ \mu = ( 4 t , 2 s ) \models 2 n $
        is of the form
        $$ P ( 2 s ) \times ( \prod _ { i \in I } ( H \wr P ( 2 ^ i ) ) $$
        where $ t $ has $ 2 $-adic expansion $ y = \sum _ { i \in I } 2 ^ i $.}
{
}
\end{block}
\setcounter{chapterno}{3} \setcounter{lastplace}{1}
\vspace*{5mm}
\begin{center} {\large {\bf \arabic{chapterno}.} {\bf The ordinary characters of the components} } \end{center}
\vspace*{-6mm}
 \begin {block}
   { So far we have
     determined the isomorphism type of the components of $ k \Xi
     {2 n} $ by identifying their vertices and Brauer quotients. We
     now turn our attention to computing the character of the complex
     lift of these components.  Let $ \Od $ be a complete
     discrete valuation ring with residue field $ k $ and field
     of fractions $ K $; then the components of $ k \Xi _ { 2 n } $
     lift to components of $ \Od \Xi _ { 2 n } $, and the
     required character of a component $ M $ of $ \Od \Xi _ { 2 n } $
     is the ordinary character of the $ K G $-module $ K \otimes_{\O} M $.
     In the following we fix $ n $ and write $ G = \Sym ( 2 n ) $.
     }  { }
 \end  {block}

     The lifts of the components of $ k \Xi _ { 2 n } $ have the same vertices
     as the mod-$2$ reductions, so we can keep the labelling:

\begin {lemma}
  {  The permutation $ { \mathcal O} \Sym ( 2 n ) $-module $ {
     \mathcal O} \Xi _ { 2 n } $ has exactly $ \lfloor n / 2 \rfloor + 1 $
     components, labelled as $ M_\mu $ where $ \mu $
     runs through the compositions of $ 2 n $ of the form
     $\mu = (4t, 2s)$ with $t, s\geq 0$. The
     vertices of these components are pairwise non-conjugate in $ G $.
     Let $ Q_\mu $ be the vertex of $ M_\mu $. Then $ Q_\mu $ is
     contained in $ \Alt ( 2 n ) $ if and only if if $ \mu = ( 4 t , 0 ) $.
     If so then $ n $ is even.
}
{
     See \arabicc{block:vertices}{}.
}
\end  {lemma}
\vspace*{0.15cm}
\newline
{{\bf Notation}
    By our labelling, the component $ M_\mu $ corresponds to the
    fixed point set $ W_\mu = ( * _ { i \in I } V_i ) * U ^ s $;
    and the vertex $ Q_\mu $
    of $ M_\mu $ is then a Sylow $ 2 $-subgroup of the stabilizer
    of $ W_\mu $. By coprimeness, there is a factorization
    $$ Q_\mu = Q_ { ( 4t , 0 ) } \times P ( 2 s ) $$
    where $ P ( 2 s) $ is a Sylow $ 2 $-subgroup of $ \Sym ( 2 s ) $
    and $ Q_ { ( 4t , 0 ) }  $ is a direct product of subgroups of the form
    $ H \wr P ( 2 ^ i ) $, as in \arabicc{block:vertices}{}.
} {}
\mylabel{lemma:alt}
\begin{lemma}
  { 
    Denote by $ \alt $ the alternating module, then 
     $$ M_{ ( 4 t , 0 ) } \otimes alt \cong M_ { ( 4 t , 0 ) } .$$ }
  { 
    Write $ \nu = ( 4 t , 0 ) $.
    Let $ X $
    be a $ Q $-invariant basis of $ M ( \nu ) $.  The vertex $ Q_{ \nu
    } $ of $ M_{ \nu } $ is contained in $ \Alt ( 4 t ) $ so
    $$
    M_{\nu } \otimes \alt \ | \Ind_Q ^ { \Sym ( 4 t ) } ( k )
    \otimes \alt \cong \Ind_Q ^ { \Sym ( 4 t ) } ( k \otimes \alt ) =
    \Ind_Q ^ { \Sym ( 4 t ) } ( k ) .$$
    Thus $ M_{ \nu } \otimes \alt $
    is a $ p $-permutation $ k \Sym ( 4 t ) $-module with vertex $ Q_{
      \nu } $.  It will suffice to show that this module has the same
    Brauer quotient (in the sense of \arabicc{block:Brauer}{}) as $ M_{ \nu } $.  Now $ Q $ only contains even
    permutations so $ \{ x \otimes 1 : x \in X \} $ is a $ Q $-invariant
    basis of $ M_{\nu } \otimes \alt $.  Thus as $ { \mathcal O} N_ { 2
      n } $-modules we have
    $$
    ( M_{ \nu } \otimes \alt ) ( Q ) \cong M_ { \nu} ( Q ) \otimes
    \alt .$$
    The images in $ { \overline N_ { \nu }} \cong \Sym ( 3 )
    \times \Sym ( 3 ) \cdots \times \Sym ( 3 ) $ of odd permutations of $
    N_{ \nu } $ are again odd permutations.  Therefore as $ { \mathcal O}
    { \overline N}_ {2n} $-modules we have
    $$
    ( M_{ \nu } \otimes \alt ) ( Q ) \cong M_ { \nu} ( Q ) \otimes
    \alt .$$
    Now $ M_ { \nu} ( Q ) = E \# E \# \cdots \# E $ so
    $$
    ( M_{ \nu } \otimes \alt ) ( Q ) \cong ( E \# E \# \cdots \# E )
    \otimes \alt \cong ( E \otimes \alt ) \# ( E \otimes \alt ) \# \cdots
    \# ( E \otimes \alt ) .$$
    But $ E \otimes \alt \cong E $, as may be
    seen by considering characters, and so $ ( M_{ \nu } \otimes \alt )
    ( Q ) \cong M_{ \nu } ( Q ) $.  }
\end {lemma}
\par
Write $ \overline M_\mu $ for the $ k G $-module $ k \otimes_{ \Od } M_\mu $.
\begin {lemma}
  { Let $ H = \Sym ( 4 t ) \times \Sym ( 2 s ) \leq G $. Then
    $ N_G ( Q_ { ( 4 t , 2 s ) } ) \leq H $ and the $ k G $-module    
    $ \overline M_{ ( 4 t, 2 s ) } $ is the Green correspondent of the
    $ k H $-module $ \overline M_{ ( 4 t, 0 ) } \# k $.  }
  {    
    A vertex $ Q $ of $ M_{ ( 4 t, 2 s ) } $ can be taken to be
    $ Q_{ ( 4 t , 2 s ) } $, recalling the notation of
    \arabicc{notation:vertex}{}. By
    \arabicc{theorm:MainTheorem}{} the set of fixed points of $ Q $
    satsifies $ X = Y * Z $ for coprime $ Y $ and $ Z $ such that
    $ Y * Z \subset H $. So
    $ N_X \leq H $ by \arabicc{lemma:permiso}{}, and we must have
    $ N_G ( Q ) \leq H $ too, since $ N_G ( Q ) \leq N_X $.

    We may therefore let the $ {
      \mathcal O} H $-module $ L $ be the Green
    correspondent of $ M_{ ( 4 t , 2 s ) } $.
    Each of $ M_{ ( 4 t , 2 s ) } , L $ and $ M_{ ( 4 t , 0 ) } \# k $ has
    vertex $ Q $, the latter by \arabicc{block:vertices}{}.
    The Brauer quotient of
    $ M_{ ( 4 t , 2 s ) } $, and hence also of $ L $, is,
    by \arabicc{block:Brauer}{},
    the $ k ( N / Q ) $-module \linebreak
    $ M_{ ( 4 t , 0 ) } ( Q ) \# k $,
    where $  M_{ ( 4 t , 0 ) } ( Q ) $ denotes the Brauer quotient
    of $ M_{ ( 4 t , 0 ) } $.
    This implies that
    $ L \cong M_ { ( 4 t , 0 ) } $ and the result follows.  }
 \end  {lemma}  
 \begin  {lemma}
   { Let $ 2 n = 4 t + 2 s $ . The $ k \Sym ( 2 n ) $-module $
     M_{ ( 4 t , 2 s ) } $ is a component of
     $ \Ind ^ { \Sym ( 2 n ) } _ H ( M_{ ( 4 t , 0 ) } \# k ) $.  }  {
     This follows immediately from the preceding lemma.  }
 \end  {lemma}  
 \begin {corollary}
   { Let $ 2 n = 4 t + 2 s $ . The module $ M_{ ( 4 t , 2 s ) } \otimes
     \alt $ is a component of \linebreak $ \Ind ^ { \Sym ( 2 n ) }
     _{ H } ( k \Xi_ { 2 n } \# \alt ) $.  }
   { By the previous lemma $ M_{ ( 4 t , 2 s )  } \otimes \alt $ is a
     component of
     $$
     ( \Ind ^ { \Sym ( 2 n ) } _{ H } ( M_{ ( 4 t , 0 ) } \# k ) ) \otimes \alt \cong \Ind ^ { \Sym ( 2 n )
     } _{ H } ( (M_{ ( 4 t , 0 ) } \otimes \alt )
     \# \alt );
     $$
     and by \arabicc{lemma:alt}{} this module is isomorphic to
     $$
     \Ind ^ { \Sym ( 2 n ) } _{ H }
     ( M_{ ( 4 t , 0 ) } \# \alt ) ,$$
     and the result is implied.  }
\end {corollary}
\begin{corollary}
  { Every irreducible character contained in the decompostion of the
    ordinary character of the module $ M_{ ( 4 t , 2 s ) } $ is indexed by a
    partition whose conjugate contains exactly $ 2 s $ odd parts.  } 
  {
    The previous corollary along with Theorem \arabicc{theorem:irs}{}
    tells us that the
    conjugate of the module $ M_{ ( 4 t , 2 s ) } $ has an ordinary character
    for which each component is indexed by a partition
    with exactly $ 2 s $ odd parts. The claim follows directly.  }
\end{corollary}
\begin{theorm}
   { Let $ 2 n = 4 t + 2 s $ . The module
     $ M_{ ( 4 t , 2 s ) } $ has character $ \phi_ {( 4 t , 2 s )} $
     given by
     
     $$
     \phi_ {( 4 t , 2 s )} = \sum_ { \lambda \in \Lambda_ { 2 n } ^ 0 \cap
       \Lambda ^{ \prime 2 s}_ {2n} } \chi ^ {  \lambda  } ,$$
     where $
     \Lambda ^{ \prime u}_ {v} $ is the set of partitions of $ v $
     whose conjugates have exactly $ u $ odd parts.  In other words,
     the module $ M_{ ( 4 t , 2 s ) } $ has character equal to the sum of
     all irreducible characters of $ \Sym ( 2 n ) $ that are
     indexed by a partition that has no odd parts but whose conjugate
     contains exactly $ 2 s $ odd parts.
     
   } { By \arabicc {theorem:summands}{}
     each module $ M_{ ( 4 r, 2 n - 4 r ) } $ occurs with
     multiplicity one as an component of $ k (\Xi_{2
       m}) $. Therefore by \arabicc{lemma:irs}{} we have the following equality of
     characters:
     $$
       \chi_ { \Xi_ {2 n} } = \sum_ { \lambda \in
       \Lambda_ { 2 n } ^ 0 } \chi ^ {  \lambda  } = \sum_ {r =
       0}^{  \lfloor n / 2 \rfloor } \phi_ { ( 4 r , 2 n - 4 r ) } .$$
     However, the preceding
     corollary tells us that $ M_{ ( 4 r , 2 n - 4 r ) } $ is a component of a
     module whose character is, by applying \arabicc{theorem:irs}{} to its
     conjugate, given by:
     $$
     \chi _ { ( 4 r  , 2 n - 4 r ) } := \sum_ { \lambda \in \Lambda ^{
     \prime ( 2 n - 4 r )}_ { 2 n } } \chi ^ { \lambda } .$$
     In other words
     there is a subset $ \Omega ^{ \prime ( 2 n - 4 r )}_ { 2 n } $ of $
     \Lambda ^{ \prime ( 2 n - 4 r )}_ { 2 n } $ such that
     $$
     \phi_ { ( 4 r, 2 n - 4 r ) } = \sum_ { \lambda \in \Omega ^{
         \prime ( 2 n - 4 r )}_ { 2 n } } \chi ^ { \lambda } $$
     so
     that
     $$
     \sum_ { \lambda \in \Lambda_ { 2 n } ^ 0 } \chi ^ { 
       \lambda  } = \sum_ {r = 0}^{ \lfloor n / 2 \rfloor }
     \phi_ { ( 4 r , 2 n - 4 r ) } =
     \sum_ {r = 0}^{\lfloor n / 2 \rfloor} \sum_ { \lambda \in \Omega ^{ \prime
         ( 2 n - 4 r ) }_ { 2 n } } \chi ^ { \lambda } .$$
     The sets $
     \Lambda ^{ \prime ( 2 n - 4 r )} _ { 2 n }
     \cap \Lambda_ { 2 n } ^ 0 $
     partition $ \Lambda_ { 2 n } ^ 0 $
     , and since every $ \Omega ^{ \prime
       ( 2 n - 4 r )}_ { 2 n } $ is a subset of $ \Lambda ^{ \prime ( 2 n - 4 r )}
     _ { 2 n } \cap  \Lambda_ { 2 n } ^ 0 $ we must have $ \Omega ^{ \prime ( 2 n - 4 r )}_ { 2 n }
     = \Lambda ^{ \prime ( 2 n - 4 r )} _ { 2 n } \cap  \Lambda_ { 2 n } ^ 0 $, completing the
     proof.}
\end {theorm}


\begin{thebibliography}{3}
  
\bibitem{bro85} M. ~Broué, \newblock On Scott modules and
  $p$-permutation modules: an approach through the Brauer morphism.
  \newblock {\em Proc. Amer. Math. Soc.}  93 (1985), no. 3, 401--408
  
\bibitem{phd04} P. Collings, \newblock Fixed point free involutions
  over a field of characteristic two and other actions of the
  symmetric group by conjugation on its own elements.  \newblock {\em
    PhD Thesis} Oxford University 2004
  
\bibitem{tp05} P. Collings, \newblock Fixed Point Sets For Permutation
  Modules \newblock {\em arXiv:0901.4378} 
  
\bibitem{gre58} J. ~A. Green, \newblock On the indecomposable
  representations of a finite group.  \newblock {\em Math. Z.  70}
  1958/59 430--445
  
\bibitem{ing90} N. F. J. Inglis, R. W. Richardson and J. Saxl
  \newblock An explicit model for the complex representations of $S\sb
  n$.  \newblock {\em Arch. Math. (Basel) 54} (1990), no. 3, 258--259

\bibitem{mur06} J. Murray, \newblock Projective Modules and Involutions.
  \newblock {\em J. Algebra} 299 (2006), no. 2, 616-622

\bibitem{the95} J.~Thévenaz,  \newblock $ G $-algebras and modular
  representation theory.  \newblock {\em {O}xford Mathematical
    Monographs, Oxford Science Publications.}  {T}he {C}larendon
  {P}ress, {O}xford {U}niversity {P}ress, New York, 1995.

\end{thebibliography}
\end{document}